\documentclass[a4paper,twoside]{article}

\usepackage{amsmath,amsthm,amsfonts,latexsym,amscd,amssymb,enumerate}

\usepackage{color}
\usepackage{color}
 \long\def\red#1{\textcolor {black}{#1}}

\swapnumbers
\theoremstyle{plain}
\newtheorem{theorem}{Theorem}[section]
\newtheorem{lemma}[theorem]{Lemma}
\newtheorem{corollary}[theorem]{Corollary}
\newtheorem{proposition}[theorem]{Proposition}
\newtheorem{conjecture}[theorem]{Conjecture}
\newtheorem{assertion}[theorem]{Assertion}

\theoremstyle{definition}
\newtheorem{definition}[theorem]{Definition}

\theoremstyle{remark}
\newtheorem{remark}[theorem]{Remark}

\newcommand{\reals}{\mathbb{R}}

\newcommand{\integers}{\mathbb{Z}}

\newcommand{\boundary}[1]{\partial#1}

\newcommand{\abs}[1]{\left\lvert#1\right\rvert} 

\newcommand{\tensor}{\otimes}

\newcommand{\iso}{\cong}

\DeclareMathOperator{\co}{co}    

\DeclareMathOperator{\supp}{supp}   

\DeclareMathOperator{\pr}{pr}

\newcommand{\forget}[1]{}

{\catcode`@=11\global\let\c@equation=\c@theorem}

\allowdisplaybreaks[2]

\newcommand{\M}{\mathcal{S}}

\newcommand{\W}{W}
\newcommand{\s}{S}
\newcommand{\X}{X}
\newcommand{\Y}{Y}
\newcommand{\Z}{Z}

\newcommand{\wE}{\widetilde{E}}
\newcommand{\Sm}{S}

\begin{document}
\pagestyle{myheadings}
\markboth{Thomas Schick, Robert S.~Simon, Stanislaw Spie\.z, Henryk Toru\'nczyk}{A parameterized version of
  the Borsuk-Ulam theorem}



\title{A parameterized version of the Borsuk-Ulam theorem}

\author{
Thomas Schick
\\
Mathematisches Institut\\Georg-August-Universit\"at\\ Bunsenstr. 3\\ 37073
G{\"o}ttingen, Germany\\ \small{schick@uni-math.gwdg.de} 
\and Robert Samuel Simon\\ Department of Mathematics\\ LSE\\
Houghton Street\\
London, WC2A 2AE,
 UK\\ \small{r.s.simon@lse.ac.uk}
\and Stanislaw Spie\.z\\ IMPAN\\ ul. Sniadeckich 8\\
P.O. Box 21\\
00-956 Warszawa,
 Poland\\ \small{spiez@impan.pl}
\and Henryk Toru\'nczyk\\ Faculty of Mathematics\\
University of Warsaw\\
Banacha 2\\
02-097 Warszawa, Poland\\
\small{torunczy@mimuw.edu.pl}
}

\maketitle

\begin{abstract} We show that for a ``continuous'' family of Borsuk-Ulam
  situations,
  parameterized by points of a compact manifold $W$, its solution set also
  depends ``continuously'' on the parameter space $W$. ``Continuity'' here
  means that the solution set supports a homology class which maps onto the
  fundamental class of $W$. When
$W\subset \reals ^{m+1}$ we also  show how to construct such a ``continuous''
family starting from a family depending in the same way continuously on the
points of $\boundary W$. This solves a problem related to a conjecture which
is relevant for the construction of equilibrium strategies in repeated
two-player games with incomplete information.

%
%

A new method (of independent interest) used in this context is a
canonical  symmetric
squaring construction in \v Cech homology with
$\integers/2$-coefficients.

\end{abstract}

\section{Introduction}
\label{sec:intro}

The classical Borsuk-Ulam theorem states that for every continuous map
$f\colon S^n\to\reals^n$ there exists $v\in S^n$ such that $f(v)=f(-v)$.  We
think of $\{v\in S^n\mid f(v)=f(-v)\}$ as the set of solutions to the
Borsuk-Ulam equation. The original Borsuk-Ulam theorem states that this set is
non-empty.

There are many generalizations of the Borsuk-Ulam theorem. In this paper, we
generalize to continuous families of Borsuk-Ulam situations and study the global
structure of the solution set (and its dependence on the parameters). At the
same time, we generalize to correspondences, i.e.~``multi-valued functions''.

More precisely, our ``family version'' of the Borsuk-Ulam equation
starts with a continuous map $F\colon W\times S^n\to \reals^n$ and
asks about properties of the solution set {$B:=\{(w,v)\in
W\times S^n\mid F(w,v)=F(w,-v)\}$}. In particular, we establish a
suitable kind of continuous dependence of the solution sets
$B_w:=\{v\in S^n\mid F(w,v)=F(w,-v)\}$ on $w\in W$. Of course,
this has to be made precise in a suitable way. Our main result,
Theorem \ref{theo:param_Borsuk_Ulam}, implies that the solution
set $B$ (as subset of $W\times S^n$) supports a (\v Cech) homology
class which is mapped to the fundamental class if $W$ is a compact
manifold. If $W=[0,1]$ this means essentially that there is a path
$s\colon [0,1]\to B$ with $s(0)\in \{0\}\times S^n$ and $s(1)\in
\{1\}\times S^n$. The homological statement is the correct
generalization if $W$ is an arbitrary manifold. \cite[Theorem 3.1]{BMZ}
is a similar family Borsuk-Ulam theorem to ours, with the precursor
\cite[Theorem 4.2]{MR1710731}. Note that the assumptions in these works are
quite different to ours, and so are the methods.

We show then that a special case of a continuous family of Borsuk-Ulam
situations arises as follows. Given is an $n$-dimensional manifold with
boundary, embedded into $\reals^n$ and on its boundary a continuous function
(or correspondence) with values in $\reals^{n-1}$. Each of its interior points
can be considered as a ``midpoint'', giving rise to a corresponding
Borsuk-Ulam equation (that the function assumes equal values on opposite points of
the boundary, ``opposite'' with respect to the
interior point). We establish in Theorem
\ref{theo:spherical fundamental class} that this is indeed a continuous family
of Borsuk-Ulam situations. As a consequence (Theorem
\ref{theo:special_fam_bors}), the solution
set has the (homological) continuity property alluded to above.

This result is strongly reminiscent to a (more complicated) Borsuk-Ulam type
statement (Conjecture \ref{conj:close}) which is part of a program in game
theory. This
program would establish the existence of equilibria in general repeated
two-player games with incomplete information. We describe the game theoretic
context in some detail in Section \ref{sec:game_appl}.

The proof of our main Theorem \ref{theo:param_Borsuk_Ulam} relies on a new
homological construction which we expect to be of independent interest,
described in Section \ref{sec:invar-homol-squar}. It is a functorial symmetric
square: a transformation $a\mapsto a^s\colon H_k(X;\integers/2)\to
H_{2k}(X^s;\integers/2)$, where $X^s=(X\times X/\tau,\Delta)$, $\tau\colon
X\times X\to X\times X$ interchanges the factors, and $\Delta$ is the image of
the diagonal. One should think of ${a}^s$ as one half of the K\"unneth
product $a\times a$ (defined on ``half of $X\times X$''). The main property we
use (and establish) states that if $X$ is a manifold and $a$ its fundamental
class, then ${a}^s$ is the fundamental class of the relative manifold
$X^s$.

The paper is organized as follows: in Section \ref{sec:param-bors-ulam-1} we
state our version of the family Borsuk-Ulam theorem and introduce the relevant
notation. In Section \ref{sec:invar-homol-squar} we discuss the symmetric
squaring transformation. In Section \ref{sec:proof1} we prove our Borsuk-Ulam
type theorem and in Section \ref{sec:game_appl} we describe the game theoretic
background and motivation.

\section{A parameterized Borsuk-Ulam theorem}
\label{sec:param-bors-ulam-1}

\label{sec:statements}

Below, $H$ denotes the {\v C}ech homology functor with $\integers /2 $
coefficients; see also Section \ref{sec:invar-homol-squar}.

\begin{definition} Let $X$ be a compact space and $(W, W_0)$ a compact
  pair. Later, $(W,W_0)$ will be a compact manifold with boundary.

  \begin{enumerate}
    \item We say that a map $f\colon X\to W$ is \emph{$H$-essential
      for}
    $(W,W_0)$ if it induces a surjection $H_d(X, f^{-1}(W_0))\to H_d(W, W_0)$
    for $d=\dim W$, compare to \cite{AH}.
  \item
    Let $p\colon E\to W$ be a fiber bundle, e.g.~$E=Y\times W$ for some
    topological space $Y$. Assume
    that $X\subset E$. We say that $X$ has \emph{property $\M$ for $(W,W_0)$}
    if the projection $p|_X\colon X\to W$ is $H$-essential for $(W,
    W_0)$. Here, $\M$ stands for ``spanning'', compare with \cite{SST1,SST2}.
 \item
     In case $W$ is a manifold with boundary $\boundary W$ (possibly empty) or $W$ is a
    compact subset of $\reals ^n$ with topological boundary $\boundary W\ne W$
    in $\reals ^n$, then in the definitions above we say ``for $W$'' in place of
    ``for $(W, \boundary W)$''.
  \end{enumerate}
\end{definition}

\begin{remark}
Suppose $W$ is a connected compact manifold with boundary $\boundary W$. The top-degree homology group of $(W, \boundary W)$
is generated by a
single element, which we call the \emph{fundamental class}
of the manifold $W$ and denote $[W]$. Thus $f$ being $H$-essential for
$(W, \boundary W)$ is equivalent to the existence of a homology class
$\alpha \in H (X, f^{-1}(\boundary W))$ such that $f_*(\alpha)=[W]$. We say
that such a class $\alpha $ \emph{witnesses} the $H$-essentiality of the map
$f$. We may also speak of an analogously defined witnessing of the property $\M$
of $X$ for $(W, \boundary W)$.
\end{remark}

\begin{definition}\label{def:Borsuk_Ulam_of_spherical}
Let  $\s=S^m $ denote the $m$--sphere. We assume that $W$ is a
compact PL--manifold  (boundary admitted), and
$\Z\subset \W\times \Sm\times \reals ^m$ is compact.

  We associate to $\Z$ its \emph{Borsuk-Ulam correspondence}
  $B(\Z)\subset W\times \reals^m $ defined by
  \begin{equation*}
    (w,e)\in B(\Z) \iff \exists x\in \s:\; (w,x,e)\in \Z\text{ and }
    (w,-x,e)\in \Z.
  \end{equation*}
\end{definition}

The first of our two results on correspondences asserts the following:

\begin{theorem}\label{theo:param_Borsuk_Ulam}
If $\Z$ has property $\M$ for
$W\times S$, then $B(\Z)$ has it for $W$.
\end{theorem}

\begin{remark}
  \begin{enumerate}
    \item Note that, if $\W=\{pt\}$ and $\Z$ is the graph of a continuous
    function $\s\to\reals^m$, then the Borsuk-Ulam theorem states that $B(\Z)$
    is 
    non-empty, whence the name chosen for the Borsuk--Ulam correspondence.
  \item
    For the same reason, Theorem \ref{theo:param_Borsuk_Ulam} is a
    parameterized Borsuk-Ulam theorem. Indeed, it yields that theorem when
    $\W=\{pt\}$, and, loosely speaking, it asserts in general the following:
    for a ``continuous'' family of Borsuk-Ulam situations, parameterized by a
    manifold
    $W$, its solution set depends continuously \red{on the
      parameters in the sense that} it
    supports a homology class which hits the fundamental class of $W$ (in
    particular, it surjects onto $W$).
  \end{enumerate}
\end{remark}

With $\pr_W$ denoting the projection to $W$ we want to describe the
construction of the class $\gamma \in H (B(\Z), \pr_W^{-1}(\boundary W)\cap B(\Z))$
satisfying $(\pr_{W})_*(\gamma)=[W]$  more precisely (which will of course
also be necessary to establish its properties). To do this, we will use  a relative
squaring construction in \v Cech homology. We believe that this construction
can be useful in other contexts and therefore deserves independent interest;
it will be described in the Section \ref{sec:invar-homol-squar}.

\begin{remark} The theorem above remains true when the manifold pair
$(W, \boundary W)$ is replaced by a \emph{relative PL--manifold}, by which we
mean a compact pair $(W, W_0)$ such that $W\setminus W_0$ is a PL-manifold.
We will skip the argument, which however will be implicit in the proof of
the result below. We get notable examples of such relative manifolds by
taking $W$ to be a compact $n$--dimensional subset of $\reals ^n,\, n<\infty$,
and $W_0$ to be its topological boundary in $\reals ^n$.
\end{remark}

\smallskip In our second theorem we assume that $W$ is  a compact subset
of $\reals^{m+1}$ of codimension zero, of course with a topological boundary in
$\reals  ^{m+1}$, which we denote $\boundary W$. We keep denoting by
$S=S^{m} $ the  $m$-sphere.

\begin{definition}\label{def:spherical_corresp}
The \emph{spherical correspondence} $\Z$ associated to a compact set
$\Y\subset \boundary W\times E$, where $E$ is an arbitrary space,
is defined by
  \begin{equation*}
      \Z:= \overline{\{(w,x,e)\in W^\circ \times \Sm\times E\mid
      (w+\lambda x,e)\in Y \text{ for some } \lambda>0\}}.
\end{equation*}
Here $W^\circ$ is the interior of the manifold $W$ with boundary, or the
topological interior of the subset of $\reals^n$, respectively. The over-line
stands, as usual, for the closure, i.e.~$\Z$ is the closure in
$W\times\Sm\times E$ of the set of points listed.

If we think of $Y$ as a multi-valued map from $\boundary W$ to $E$, $Z$ is a
canonical ``extension'' to $W$ using the convex structure of $\reals^{m+1}$
and keeping track of the directions needed for the extension. It collects all
Borsuk-Ulam equations arising if the points in $W$ are treated as ``midpoints''
of $W$.
\end{definition}

\begin{theorem}\label{theo:spherical fundamental class}
  If $\Y\subset \boundary W\times E$ has property $\M$ for
  $(\boundary W, \emptyset) $ then  the spherical correspondence
  $\Z\subset W\times S\times E$ associated to $\Y$ as above has
  property $\M$ for
  $W\times S$.
\end{theorem}

As a consequence, we immediately obtain from Theorem \ref{theo:spherical
  fundamental class} and Theorem \ref{theo:param_Borsuk_Ulam} \red{the following
corollary.}

\begin{corollary}\label{theo:special_fam_bors}
  Let $W\subset\reals^n$ be a compact $n$-dimensional embedded manifold with
  boundary and assume that $Y\subset \boundary W\times \reals^{n-1}$ has
  property $\M$ for $(\boundary W,\emptyset)$.

  Then $Z:=\{(x,v)\in W\times\reals^{n-1}\mid \exists x_1,x_2\in\boundary W:\,
  x\in [x_1,x_2],\,(x_1,v),(x_2,v)\in Y\}$ satisfies property $\M$ for
  $(W,\boundary W)$.
\end{corollary}

For earlier results related to {Corollary}
\ref{theo:special_fam_bors} see \cite{Joshi,Ol}.

\section{\v Cech homology and symmetric homology squaring}
\label{sec:invar-homol-squar}


Throughout this note, all manifolds are  finite-dimensional, possibly with a
nonempty boundary, and all spaces encountered will be subspaces of
manifolds. The homology groups we are using will exclusively be
\emph{\v Cech homology} groups. Their properties can be found in
\cite[Chapters IX, X]{MR14:398b} and in \cite[VIII,13]{MR96c:55001}. We list the most important properties (not all of them will
be relevant to us):
\begin{enumerate}
\item \v Cech homology is defined for
compact pairs.

\item \v Cech homology satisfies excision in a very strong form: if
  $f\colon (X,A)\to (Y,B)$ is a map of compact pairs such that $f|\colon
  X\setminus A\to Y\setminus B$ is a homeomorphism, then $f_*\colon
  H (X,A)\to H (Y,B)$ is an isomorphism of \v Cech homology groups.
\item More generally, the \emph{Vietoris theorem} about maps with acyclic
  fibers holds: if $f\colon  (X,A:=f^{-1}(B))\to (Y,B)$ is a surjective map of compact pairs such
  that the reduced \v Cech homology groups $\tilde H (f^{-1}(y))$
  are trivial for all $y\in Y\setminus B$, then $f_*\colon H (X,A)\to H (Y,B)$
  is an isomorphism of \v Cech homology groups \cite{Vietoris} (using
  \cite[Theorem 5.4]{MR14:398b}).
\item For a Euclidean neighborhood retract (ENR), e.g.~for a topological
  manifold, \v Cech homology and 
  singular homology are canonically isomorphic.
\item \label{item:continuity}Continuity: The \v Cech homology functor commutes with
taking an inverse limit.
\item \label{item:intersection pairing} For \v Cech homology, there is a
  natural intersection pairing: If  $(X,A)$ and $(Y,B)$ are two compact pairs
  inside an $m$-dimensional manifold $W$ (possibly with non-empty boundary),
then there is  an intersection pairing
  \begin{equation*}
H_p(X,A)\tensor H_q(Y,B)\ni \alpha \tensor \beta\ \mapsto \alpha \bullet \beta \in
H_{p+q-m}
((X,A)\cap (Y,B))
  \end{equation*}
  which is natural for inclusions of pairs.
(Here, $(X,A)\cap (Y,B):=(X\cap Y, A\cap Y \cup X\cap B)$.)
The reasoning  in \cite[p. 342]{MR96c:55001} not  only shows that the product  $\alpha\bullet \beta $ can be computed in any
  given neighborhood of $X\cap Y$, but actually that it can be defined whenever
  a neighborhood of $X\cap Y$ is a manifold (even if $W$ isn't).

\item For \emph{compact} subsets of manifolds, and with coefficients
  in a field (e.g.~$\integers/2$), \v Cech homology is a homology theory,
  in particular with a long exact sequence of a pair.
  This does not  hold in general for $\integers$-coefficients.
  The main reason however for which we work with $\integers/2$-coefficients
  is to avoid orientability assumptions.
\end{enumerate}

As in \cite{SST2} these properties allow to define a certain restriction operator, as
follows. Suppose $(X, A)$ and $(Y, B)$ are compact pairs in
a space $Z$. If $X\cap (Y\setminus B)$ is a relatively open subset of
$X\setminus A$ then we say that $(Y,B)$ is \emph{admissible} for $(X,A)$.
In this case we have canonical homomorphisms
\[H (X, A)\to H (X, X\setminus (Y\setminus B))\to H (X\cap Y, X\cap B)\] the
first of
which is induced by the inclusion and the second is the excision of
$X\setminus Y$ from $X$.
The arising composition $H (X,A)\to H (X\cap Y,X\cap B)$ will be called
\emph{restriction} and denoted $\alpha\mapsto \alpha{|(Y,B)}$.
Note that then any other compact pair $(Y',B')$ in $Z$ such that
$X\cap (Y',B')=X\cap (Y,B)$ is admissible for $(X,A)$, too,  and
$\alpha {|(Y',B')}=\alpha{|(Y,B)}$.

\bigskip From now on we assume that $H$ denotes the graded \v Cech homology
functor with $\integers/2$ coefficients. Our aim is to construct a ``symmetric squaring'',
similar to the squaring $H_k(X,A)\to H_{2k}(X\times X, A\times X\cup X\times A)$
obtained from the exterior  homology product \cite{MR96c:55001}, but which
takes values in the homology of the (slightly modified) {\emph{symmetric}} square
of the pair $(X,A)$.  We explain this below.

\smallskip Let $\X\supset A$ be compacta. On $\X\times \X$ we have the
coordinate-switching involution $\tau$. We put
 {$$(\X, A)^s=\Big(\X\times \X/ \tau, \big(\Delta \cup (A\times \X)\cup (\X\times A)\big)/\tau\Big)
 \, \text{ and }\, \X^s=(\X, \emptyset )^s $$}
where $\Delta =\{(x,x)\colon x\in X\}$ is the diagonal. Clearly, a map
$f\colon (\X,A)\to (\Y, B)$ of compact pairs induces a map $(\X,A)^s\to (\Y,B)^s$,
denoted $f^s$. (Here and below the superscript $s$ stands for ``symmetric square''.)

\begin{theorem}\label{theo:quotientproduct}For each $k$ there is an 
assignment  $$H_k(\X,A)\ni\alpha\mapsto \alpha ^s\in H_{2k}((\X, A)^s)
$$
such that
   \begin{enumerate}

\item\label{item: naturality} It is natural: if $f\colon (\X, A)\to (\Y, B)$
is a map of compact  pairs then
$$(f_*(\alpha))^s=f^s_*(\alpha ^s)\text{ for } \alpha\in
H_k(\X,A)$$

  \item If $(\X,A)=(M, \boundary M)$ for some compact PL--manifold $M$ and
  $\alpha =[M]$ is its fundamental class, then $\alpha ^s$ is the
  fundamental class of the relative PL--manifold $(M, \boundary M)^s$.
   \end{enumerate}
\end{theorem}

Before proving this 
theorem
we need also to explain what we mean by a
\emph{fundamental class} of a relative PL-manifold. When $Y$ is a genuine
connected manifold and $B=\boundary Y$ is its boundary (possibly empty), then
the fundamental class $[Y]\in H_m(Y,B),\,m=\dim Y$, has already been defined
above. It is unique,
because we use $\integers/2$ as the coefficients.  In case $Y$ is disconnected
but remains
an $m$-manifold we define $[Y]\in H_m(Y, \boundary Y)$ as the sum of the
images of the classes $[Y_i]$ under the homomorphisms induced by inclusions
$(Y_i,\boundary Y_i)\hookrightarrow (Y, \boundary Y)$, where $Y_i$ runs over
the (finitely many) components of $Y$.
It is clear that when $M$ is a codimension 0 submanifold of $Y$,
contained in the interior of $Y$, then $[M]$ is the restriction of $[Y]$
to $(M, \boundary M)$.
Thus in the general case  of a relative PL-manifold $(Y,B)$ we may exhaust
$Y\setminus B$ by compact $m$--manifolds $M_i, i\in \integers$, each
contained in the interior of the next, and define the fundamental class of $(Y,B)$
as $\alpha =\lim\limits_{\longleftarrow} \alpha _i\in H_m
(Y,B)$,
where
$\alpha _i\in H_m(Y, Y\setminus M_i^\circ )$  is the image of $[M_i]$ under
the map 
$H_m(M_i,\boundary M_i) \to H_m(Y, Y\setminus M_i^\circ )$ induced
by the inclusion $(M_i,\boundary M_i) \hookrightarrow (Y, Y\setminus M_i^\circ )$.
This fundamental class $\alpha$ is uniquely characterized by the property that
$\alpha {|(M,\boundary M)}\ne 0$,
for each compact $m$--manifold $M\subset Y\setminus B$.

A word should be said about the exhaustion of $Y\setminus A$ by the
manifolds $M_i$. In the presence of a PL-structure on $Y\setminus A$ its
existence is obvious. However, for $m\ge 6$ it exists also when $Y\setminus A$
is a topological manifold, see \cite[p.
108
]{KS}.


\medskip\emph{Proof of  Theorem {\ref{theo:quotientproduct}}}. We assume first
that $X$ and $A$ are compact polyhedra, in which case \v Cech and
singular homology of $(X,A)$ coincide. Let $U$ be a neighborhood
 of  the diagonal $\Delta $ in {$X\times X / \tau$} and let
 $$(X,A)_U^s=(X,A)^s\cup ({X\times X / \tau },U).$$ 
Let $\sigma =\sigma_U$ be a
representative of $\alpha\in H_k(X,A)$ by a singular chain
$\sigma=\sum_{i=1}^n \sigma_i$ such that for all $i,j$ the image
of $\sigma_i\times \sigma_j$  under the projection $p\colon
 X\times X\to {X\times X / \tau }$ is either contained entirely in $U$ or does not
 meet the diagonal $\Delta$ at all. Then we define $\sigma ^s $ to
be  the chain $\sum_{i<j} p_*(\sigma_i\times \sigma_j)$ in  the
relative chains of $(X,A)_U^s$.  Technically speaking, we have to
subdivide $\sigma_i\times \sigma_j$ into simplices. This can be
done in an arbitrary fashion ---since we are dealing with
$\integers/2$-coefficients, not even the
  orientations play a role. Note that
  $p_*(\sigma_i\times \sigma_j)=p_*(\sigma_j\times \sigma_i)$
  since we quotient out by $\tau$.  Thus
  $$
    \boundary (\sigma ^s )=
    \sum_{i< j} p_*(\boundary \sigma_i\times \sigma_j)+
  p_*(\sigma _i \times \boundary \sigma_j)=
\sum _i p_*(\sigma _i\times \sum _{j\ne i}\boundary \sigma _j)$$

By the smallness condition imposed on the $\sigma _j$'s each
 ${p_*(\sigma _i\times \boundary \sigma _i}
)$
lies in $U$. Hence
$\boundary (\sigma ^s )=\sum _{i}p_*(\sigma _i\times  \sum _j \boundary \sigma _j)=0$
as a chain in $(X,A)^s_U$, since
$\sum _i\sigma _i $ is a cycle. This shows that $\sigma ^s _U$ is a cycle
in $(X,A)^s_U$.

One checks immediately that its homology class (actually, even the cycle
  itself) does not depend on the ordering of the simplices $\sigma_i$ chosen
  above.

  Similarly, when $\mu=\mu_U$ is a second ``$U$--small'' representative of
  $\alpha $ then   $\sigma ^s _U$ and $\mu ^s _U$ turn out to represent
  the same element of $H_{2k} ((X,A)^s_U)$.
  Indeed, we may take a $(k+1)$-chain $\sum {\tau_k}$ such that
  $\mu =\boundary (\sum_n \tau_n) + \sigma$ and assume after subdivision
  that   $\tau $ satisfy the same   ``smallness'' condition as $\sigma $ and $\mu $.
  Clearly it suffices to treat the case where there is just one simplex $\tau$, and
  to argue by induction. Assuming that in the representation
  $\mu=\sum _i\sigma _i+\boundary \tau$ all the $\sigma _i$'s
  precede the simplices $\alpha _k$ of $\boundary \tau $ we get $\mu ^s=
  \sigma ^s+ p_*(\sigma \times \boundary \tau )+(\boundary \tau) ^s$, where all the
  summands  $p_*(\alpha _k\times \alpha _l) $ of $(\boundary \tau )^s$
  lie in $U$. Since
  moreover $p_*(\sigma \times \boundary \tau )=
  \boundary (p_*(\sigma \times \tau ))$ (because $\boundary \sigma =0$) it
  follows that $\mu ^s$ and $\sigma ^s$ are homologous in
  $(X,A)^s_U$.

  Now, the classes $[\sigma _U]^s$ form an element of the inverse system
  $\{H_{2k}((X,A)^s_U)\}$, with homomorphisms induced by the inclusions.
  (The $U$'s run over all the neighborhoods of the  diagonal.)
  Since $\lim\limits_{\longleftarrow} H((X,A)^s_U)=H((X,A)^s)$ by (6),
  this system uniquely defines an element $\alpha ^s\in H_{2k}((X,A)^s)$.

This concludes   the construction in case of polyhedral pairs $(X,A)$ and it
is clear that in case $X$ is a PL--manifold  and $\alpha$ is its
fundamental class then $\alpha ^s$ is by construction the fundamental
class of $(X,\boundary X)^s$.

When $(X, A)$ is an arbitrary compact pair, 
we may represent it as an inverse limit of compact polyhedral pairs $(X_i, A_i)$.
Given $\alpha \in H_k(X,A)$ let us denote its image in $H_k(X_i,A
_
i)$ by $\alpha _i$. Then $(\alpha _i^s)$ is an element of
$\lim\limits_{\longleftarrow}
H_{2k}((X_i,A_i)^s)$,
whose limit is the class $\alpha ^s\in H_{2k}((X,A)^s)$ we are after.  By
  the very construction,  naturality is true even on the level of cycles.
\hfill    $\square$


\begin{remark} Above, one can replace the field $\integers / 2$ by other
coefficients, over which suitable manifolds are orientable.
  More details about ``invariant homology squaring'' and about the parameterized
  Borsuk-Ulam theorem can be found in the G\"ottinngen diploma thesis of
 Denise Krempasky (maiden name Nakiboglu)
 \cite{nakiboglu07:_dachab_cech_homol}. Moreover,
  \cite{nakiboglu07:_dachab_cech_homol} contains a generalization of the
  construction to homology with integer coefficients in even dimensions.
%
 It is an interesting problem to lift the invariant squaring
  construction to the generalized homology theories ``non-oriented bordism''
  or even to oriented bordism. This is work in
  progress \cite{Krempanski}.
\end{remark}

Finally, we need to establish that restriction is compatible with some other
operations.

\begin{lemma}\label{lem:naturality of restriction}
  \begin{enumerate}\item \label{item:a}
    Let $(X, A)$ and $(X',A')$ be compact pairs in spaces $Z$ and
    $Z'$, respectively, and let a map $f: Z'\to Z$ satisfy $f(X')\subset X$
    and $f(A')\subset A$.  Then for any compact pair $(Y, B)$ in $Z$,
    admissible for $(X,A)$, one has
    \[g_*(\alpha){|(Y,B)}=h_*(\alpha{{|(f^{-1}(Y),f^{-1}(B)))}} \text{ for }
    \alpha \in H(X',A')\] where $g\colon (X', A')\to (X,A)$ and $h\colon
    {X'\cap(f^{-1}(Y),f^{-1}(B))\to X\cap}(Y,B)$ are the restrictions of $f$.
  \item \label{item:b}
    Let $(X, A)$ and $(Y, B)$ be compact pairs in a space $Z$. If $(Y,B)$
    is admissible for $(X,A)$, then $(Y,B)^s$ is admissible for $(X, A)^s$ and
    $(\alpha{|(Y,B)})^s=\alpha ^s{|(Y,B)^s}$ {for $\alpha\in H(X,A)$}.
  \item \label{item:c}
     Let $(X,A), (X',A')$ and $(Y,B)$ be compact pairs in a relative
    manifold $(Z, Z_0)$. If $X\cap X'\cap Z_0=\emptyset$ and $(Y,B)$ is
    admissible for both pairs $(X,A)$ and $(X_0,A_0):=(X, A)\cap (X',A')$,
    then $(\alpha'\bullet \alpha)|{(Y,B)}=\alpha'\bullet (\alpha{|(Y,B)})$ for
    $\alpha\in H(X,A)$ and $\alpha '\in H(X',A')$.
  \end{enumerate}
\end{lemma}

\begin{proof} The proofs of the three parts are similar. For \ref{item:a},
the conclusion follows from the definition of the restriction
and the commutativity of the diagram induced on the homology level by
the following  one:
\begin{equation*}
  \begin{CD}
    (X',A') @ >>> (X',X'\setminus (f^{-1}(Y)\setminus f^{-1}(B))) @<<<
    {X'\cap (f^{-1}(Y),}f^{-1}(B)) \\
   @VVV  @VVV @ VVV \\
    (X,A) @ >>> (X,X\setminus(Y\setminus B)) @<<< {X\cap( Y,} B)
  \end{CD}
\end{equation*}
where horizontal maps are  inclusions and the vertical ones are induced by $f$.

For \ref{item:b}, we  consider a similar diagram
\begin{equation*}
  \begin{CD}
    H(X,A) @ >>> H(X,X\setminus(Y\setminus B)) @<<< H({X\cap( Y, B)}) \\
   @VVV  @VVV @ VVV \\
    H((X,A)^s) @>>> H((X,X\setminus(Y\setminus B))^s) @<<< H(({X\cap( Y, B)})^s)
  \end{CD}
\end{equation*}
where horizontal arrows are inclusion-induced and vertical ones are
given by the assignment $\alpha \mapsto \alpha^s$.   By \ref{item: naturality}
in Theorem \ref{theo:quotientproduct} the diagram is commutative, whence
again the conclusion follows from the definition of the restriction.

Finally, \ref{item:c} follows in the same manner using the commutative diagram
\begin{equation*}
  \begin{CD}
    H(X,A) @ >>> H(X,X\setminus(Y\setminus B)) @<{\iso}<< H({X\cap( Y, B)}) \\
   @VVV  @VVV @ VVV \\
    H(X_0,A_0) @ >>> H(X_0,X_0\setminus(Y\setminus B))
    @<{\iso}<< H({X_0\cap( Y, B)})
  \end{CD}
\end{equation*}
where the vertical arrows are $\bullet$--multiplications by a fixed
$\alpha '\in H(X',A')$. \end{proof}

\section{Proof of the main results}

\subsection{Proof of Theorem \ref{theo:param_Borsuk_Ulam}; the
naturality of the witnessing class.}
\label{sec:proof1}

We use the  notation of Section \ref{sec:param-bors-ulam-1}
concerning Theorem \ref{theo:param_Borsuk_Ulam} and denote
$\reals ^m$ by $E$. We consider $E$ as being embedded in the
 sphere $\wE=E\cup \{\infty\}$.
 Hence $Z\subset W\times S\times
\wE$ and $Z^s\subset (W\times S\times \wE)^s$. The last space also
contains the ``anti-diagonal'' $A=\{[(w,x,e),(w,-x,e)]: w\in W,
x\in S, e\in \wE\}$, with well-defined projections $A\to W\times
\wE$ and $p_W\colon A\to W$. Note that there is an obvious
 { mapping} $A\cap Z^s\to B(Z)$ which preserves the projections to
 $W$. {(Here and below we identify
 the pair $A\cap Z^s=(A\cap (Z\times Z/\tau), \emptyset)$ with the
 space $A\cap (Z\times Z/\tau )$.)}
 Hence it suffices to show that the projection $p_W {|A\cap
Z^s}$ is $H$-essential, i.e.~sends some homology class $\beta \in
H_k ({(A, p_W^{-1}(\boundary W)})\cap Z^s)$ to $[W]$.

To this end let $k=\dim W$ and $m=\dim S=\dim E$ and choose
$\alpha \in H_{k+m}(Z,  (\boundary W\times S\times \wE)\cap Z)$
which projects to $[W\times S]\in H_{k+m}((W, \boundary W)\times
S)$. Note that $A$ is over $W\times \wE$ homeomorphic to the
$(2m+k)$--manifold $W\times P\times \wE$, where  $P$ is the real
projective $m$-space, and $(A,\boundary A)\subset ((W,\boundary
W)\times S\times \wE)^s$. We  claim that one can take
$\beta:=[A]\bullet \alpha ^s$, where $[A]\in H_{2m+k}(A, \boundary
A)$ is a fundamental class and $\bullet $ is taken in the relative
$2(2m+k)$--manifold
$(W\times S\times \wE)^s$. Note that the product 
$[A]\bullet \alpha ^s$
makes sense,
because the anti-diagonal $A$ misses the singular set of this relative manifold.

Indeed, $Z$ is a compact set in $W\times S\times E$ and{, by
its choice, $\alpha$ is equal to
 $[W\times S\times \{0\}]$ in $H_{k+m} ((W, \boundary W)\times S\times \wE)$, by
which we mean that the images of these homology classes under the
inclusion-induced homomorphisms coincide.}
 Similarly, $[W\times S\times \{0\}]$ is
 in $H_{k+m} ((W, \boundary W)\times S\times \wE)$ equal to the fundamental
 class of the manifold $M:=W\times \{(x,u(x)\}_{x\in
S}$, where for $u\colon S\subset \reals^{m+1}\to E=\reals^m$  we
choose the projection onto the first $m$ coordinates. (The
advantage of this $u$  will be seen below.) Thus $\alpha
^s=[M{]^s}$ in $H_{2(k+m)}\left(\big{(}(W,\boundary W)\times
S\times \wE\big{)}^s \right)$ and $[A]\bullet \alpha ^s=[A]\bullet
[M{ ]^s }$ in $H_j(A,\boundary A)$, with
$j=2(k+m)+(2m+k)-2(2m+k)=k$. But{, by Theorem
\ref{theo:quotientproduct}, $[M]^s$ is the fundamental class
 of the relative manifold $(M, \boundary M)^s$. The latter intersects
 $(A, \boundary A)$ transversally along $(N, \boundary  N)$, where
  the manifold $N:=\{[(w,e_{m+1},0),(w,-e_{m+1},0)]\colon w\in W\}$
is disjoint from the singular set $\Delta $ of $M\times M/\tau $.}
 Hence $\beta =[A]\bullet
[M{]^s}=[N]$ in $H_k(A, \boundary A)$ and $(p_W)_*(\beta
)=(p_W)_*([N])$ in $H_k(W, \boundary W)$. Our claim follows, for
$p_W$ maps $N$ homeomorphically onto $W$.
\hfill
   $\square$

\smallskip\label{rem:witness formua}
\begin{remark} In fact we established the following: if $\alpha $ witnesses
property $\M$ of $Z$ for $W\times S$, then  an explicit class
$\gamma _\alpha $ witnessing property $\M$ of $B(Z)$ for $W$ may
be defined by $\gamma _\alpha =p_* (\beta)=p_*([A]\bullet \alpha
^s)$, where $p\colon {(A, p_W^{-1}(\boundary W))
\cap Z^s}\to (B(Z), (\boundary W\times E)\cap B(Z))$ is induced by
the projection $q\colon A\to W\times \wE$.
\end{remark}

\smallskip Below we show that the above dependence between the witnesses
$\gamma _\alpha$ and $\alpha $ satisfies a certain naturality property.

\begin{proposition}
If $V$ is a top-dimensional compact submanifold of $W$ then the class
$\gamma _{\alpha{|(V,\boundary V)\times S\times \wE}}$ is the restriction of
$\gamma _\alpha$ to $(V,\boundary V)\times \wE$.
\end{proposition}
\begin{proof}
Implicit in the formulation above is that $\alpha{|(V,\boundary V)\times S\times \wE}$
is a witness of the property $\M$ of $Z\cap p_W^{-1}(V)$ for $V$; compare to
\cite[Lemma 2a)]{SST2}. The latter follows also from part \ref{item:a} of Lemma
\ref{lem:naturality of restriction}, for $[V]$ is the restriction of $[W]$ to
${(V, \boundary V)}$.

For the same reason, to prove the desired equality it suffices to show that
$$([A]\bullet \alpha ^s)|(A_V, \boundary A_V)=[A_V]\bullet
(\alpha{|(V,\boundary V)\times S\times \wE})^s$$
where we write $A_V$ for $q ^{-1}(V\times \wE)$.

But as noted after the definition of the restriction in Section
\ref{sec:param-bors-ulam-1},
the pair $(A_V, \boundary A_V)$ may be replaced by
$\big((V, \boundary V)\times S\times \wE\big)^s$.
By  parts \ref{item:b} and \ref{item:c} of
Lemma \ref{lem:naturality of restriction} the left hand side is thus equal to
$([A]|(A_V, \boundary A_V))\bullet (\alpha|(V, \boundary V)\times S\times
\wE)^s$.
The result follows, for $[A]|(A_V, \boundary A_V)=[A_V]$.
\end{proof}

\subsection{Transport of the spanning property: proof of Theorem
\ref{theo:spherical fundamental class}}
\label{sec:transp-spann-prop}

  We now use the notation of Theorem \ref{theo:spherical fundamental class}.
  For the proof of this theorem we consider a number of intermediate
  correspondences. First, by the K{\"u}nneth formula, it follows that
  the correspondence $Y\times W\subset \boundary W\times E\times W$
  has property $\M$ for $\boundary W \times W$.

Next, we define the correspondence
$\widetilde{Z}\subset
W\times \Sm \times Y $
by
\begin{equation*}
    \widetilde{Z}:=
    \overline{\{(v,x,(w,e))\in W\times \Sm \times Y :
     v\ne w, x={(w-v)/{\abs{w-v}}}\}}.
\end{equation*}
and set $Z_0:=
\widetilde{Z} \cap (\boundary W\times \Sm \times Y)$.
The map
$\widetilde{Z} \to W \times Y$,
induced by the projection,
restricts to a homeomorphism
$\widetilde{Z}\setminus Z_0 \to
(W\setminus\boundary W)\times Y$
and therefore induces a {\v C}ech homology isomorphism
\begin{equation*}
  H (\widetilde{Z}, Z_0)\to
  H ((W,\boundary W)\times Y).
\end{equation*}

Define $V\subset W\times \Sm \times \boundary W$ as
\begin{equation*}
  V:=\overline{ \{(v,x,w)\in W\times \Sm \times \boundary W:
  v\ne w, x= {(w-v)/{\abs{w-v}}} \}}
\end{equation*}
and set
$V_0:= V \cap (\boundary W\times \Sm \times \boundary W)$.
As above, the projection $$ p\colon V\to W\times \boundary W$$
restricts to a homeomorphism
$V\setminus V_0 \to
(W\setminus\boundary W)\times \boundary W$
and therefore
induces a {\v C}ech homology isomorphism
$H (V,V_0)\to H ((W,\boundary W)\times \boundary W).$

Now consider the commutative diagram
\begin{equation*}
  \begin{CD}
    (W,\boundary W)\times Y @>>> (W,\boundary W)\times \boundary W\\
    @AAA @AA p A\\
    (\widetilde{Z},Z_0) @>>> (V,V_0)\\
    @VVV @VV q V \\
    (Z,\boundary Z) @>>> (W,\boundary W)\times \Sm .
  \end{CD}
\end{equation*}
where the lower vertical maps are induced by the projection
along $\boundary W$. We already know that the upper
 vertical maps induce {\v C}ech homology isomorphisms and that
the upper horizontal map induces an epimorphism on $H_{2m+1}$.
Hence so does the middle horizontal map and to prove this for the
lowest one it remains to show the following

\smallskip \noindent{\bf Claim.} {\it With the above notation, $H_{2m+1}(q)$
is an epimorphism.}

\begin{proof}
Without loss of generality we can assume that $W^\circ $ is connected.

 Fix $v$ in the interior of $W$ and let $D_v$
  be a small disc centered at $v$ contained in the interior
  of $W$, and with boundary the sphere $S_v$.
  Consider the 
  commutative
  diagram
\begin{equation*}
  \begin{CD}
  H_{2m+1}  \left( (W,\boundary W)\times \boundary W\right) @>>>
   H_{2m+1}  \left( (D_v,S_v)\times \boundary W\right)\\
    @AAA
    @AAA\\
   H_{2m+1}  \left( V,V_0 \right)
@>>>
    H_{2m+1}  \left( p^{-1}(( D_v, S_v)\times \boundary W)\right)\\
   @VV q_* V
   @VV f_* V\\
   H_{2m+1}  \left( (W,\boundary W)\times \Sm\right)  @>>>
   H_{2m+1}  \left( (D_v,S_v)\times \Sm\right)
  \end{CD}
\end{equation*}
where the horizontal homomorphisms are restrictions, the upper vertical
ones are induced by $p$, and $f$ is the restriction of  $q$.
The upper vertical arrows are isomorphisms. So are
the horizontal ones, the reason for the upper and lower ones being that
restriction takes the generator $[W]$ of
$H_{m+1}(W, \boundary W)$ to a generator of $H_{m+1}(D_v, S_v)$.

\smallskip Therefore, to prove the claim it suffices to show that the map
\begin{equation*}
g
\colon (D_v,S_v)\times \boundary W \to
 (D_v,S_v)\times \Sm ;\ (u,w)\mapsto (u,{\frac{w-u}{\abs{w-u}}})
\end{equation*}
induces an epimorphism on $H_{2m+1}$.  
(We are using the fact that $f$ is the composition of $g$ with $p$ restricted
to the set $p^{-1}(D_v\times \boundary W)$, on which it is a homeomorphism.)
However,  $f$ is homotopic to the map
\begin{equation*}
 g' \colon (D_v,S_v)\times \boundary W \to
 (D_v,S_v)\times \Sm ;\ (u,w)\mapsto (u,{\frac{w-v}{\abs{w-v}}}).
\end{equation*}
By Borsuk's separation criterion \cite[p.302]{MR14:398b} 
the map $\boundary W \to \Sm $, $w \mapsto
{\frac{w-v}{\abs{w-v}}}$, induces an epimorphism in {\v C}ech
homology $H_{m}$. Thus, by the K\"{u}nneth formula, $H_{2m+1}(g')$
is an epimorphism. Hence the same follows for $H_{2m+1}(g)$, as
desired.
\end{proof}

\section{\red{Relations} to game theory}\label{sec:game_appl}



 The inspiration for this paper came originally from
 games of incomplete information.  Infinitely repeated
             two-person, non-zero-sum games of
              incomplete information on one side were introduced
 by Aumann, Maschler, Stearns \cite{aumann}.
There is
 a finite set $K$ of states of
   nature and two players.  Nature chooses a
   state $k\in K$
 according to a commonly known probability distribution $p_0$ on $K$.
    The first player, but not the second player,
   is informed of nature's choice. The finite sets of moves for the players
 are the same for all states.
 The chosen state remains constant
  throughout the play.
Although the chosen state $k$, along with the moves of the players, determines
 the stage payoffs,
    during the play the second player learns nothing about
      his payoff, as this could give him information about
       the state of nature.

 The first player and second player
   have the  finite sets  $I$ and $J$ of moves, respectively.

  For any finite set $S$ let $\Delta (S)$ be the  probability
  simplex {$\{ \lambda \in [0,1]^S\ | \ \sum_{s\in S} \lambda^s
  =1\}$.}


 We assume that there
 are two sets of finite signals $R$ and $S$, received by the
       first and second players, respectively, and
 a stochastic signaling function,
     $\Lambda\colon K\times I\times J\rightarrow \Delta (R\times S)$.
     After each stage in
      which $i\in I$ and $j\in J$ were played, and $k$ is the state of nature,
      a member of $R$ and a member of $S$ is determined by
     $\Lambda (k,i,j)$ and communicated to Player One
     and Player Two, respectively.
     The only knowledge the players have of the moves of their opponents
      is through their observations of the sets $R$ and $S$, and the only
       knowledge Player Two has of the state of nature is from the initial
     probability $p_0$ and the received sequence of signals in $S$.

An equilibrium of the  game is a pair of strategies such that
for every state $k$ there are limits $a^k$ and $b^k$ (as the number $n$
 of stages
 goes to infinity) for  the
   accumulated averages  up to the stage $n$
 of  the expected
 payoffs of Players One and Two, respectively,
 and neither player can obtain a higher limit superior of this same
  average payoff
  (for Player Two as  weighted by  the initial probability
 distribution on $K$) by choosing a different strategy.

 \begin{conjecture}\label{conj:main}
   For all signalling functions $\Lambda$ the corresponding infinitely
   repeated game of incomplete information on one side has an equilibrium.
 \end{conjecture}

Obviously, this is the principal conjecture in this field. In full generality,
it is still open. Special cases of increasing generality have been established
by Simon, Spie\.z, and Toru\'nczyk \cite{SST1}, Renault \cite{Renault} and
again Simon, Spie\.z, and Toru\'nczyk \cite{SST2}. \red{By a suitable construction,
 to the most general game considered in the present paper one can assign a game
 as treated by \cite{SST2} with an essential
  preservation of  the payoff structure. Therefore, for those
only interested in economic applications the \red{full generality of
Conjecture \ref{conj:main}} will be unnecessary.}

\smallskip
How do  equilibria in the non-zero-sum games
  relate to
topological theorems resembling the  Borsuk-Ulam Theorem?

The general approach is the construction of rather simple strategies,
depending on a small set of parameters. Then, conditions for such a strategy
to be an equilibrium are derived. Finally, one has to show that among the
strategies considered  there is at least one which satisfies these conditions.
Such constructions have already been considered in
\cite{aumann,Renault,simon03:_habil,SST1,SST2}.

Given that the two players know each other's strategies \red{(an assumption
that can be made of an equilibrium)},
 the critical information
 used by  the second player in the special strategies considered is the
 conditional probability on the
 set $K$ of states that is induced by the long term observable  behavior of the
 first player. This conditional probability defines a multi-dimensional
martingale (as many dimensions as are states in $K$) starting at
 the initial probability $p_0$. We restrict ourselves to strategies where this
martingale is controlled entirely by the first player and that it
 converges after only a few initial stages  to some point in a
 finite  subset $V\subset \Delta(K)$ of probabilities for the state of
 nature. Two conditions together imply that the (constructed) strategy is an
 equilibrium. The first, ``individual rationality'' for each $v\in V$, is
 classically under control \cite{aumann} and we don't discuss it here.

 The second condition
 is called ``incentive compatibility''. Incentive compatibility
requires that  the first player has no
 incentive to infer  through her behavior  one  of the $v\in V$  over
 another member of $V$. If the vector $y\in {\reals}^k$ represents the expected
 payoff to the first player starting from the initial probability $p_0$ and
 $y_v$ the expected payoff to the first player after
  $v$ is chosen, this means that $y_v^k =y^k$ \red{if $v^k>0$} (with some modification if
  $v^k=0$ which \red{leads to the introduction of slack variables and the
  requirement of saturation below}). Hence we are describing
  a set $V$ whose convex hull includes
 the initial probability $p_0$ with corresponding vectors
$(y_v \ | \ v\in V)$ that are all equal. With some form of continuity of these
payoff vectors we recognize a relation to the
 Borsuk-Ulam Theorem. \red{This relation becomes stronger if there are restrictions on how we can
convexify members of V, for example that the initial $p_0$ must be in the
convex combination of at most two members of $V$.}

 Continuing with the game context,  a family ${\cal L}$ of subsets
of $K$ is introduced. A member $L$ of ${\cal L}$ is
 a set of states for which  the first player can confirm that the chosen
  state belongs to $L$
 without revealing anything more about the chosen state. A choice of an $L$ in
 $\cal L$ (which must contain the chosen state) by the first player
 determines a kind of sub-game, which motivates the mathematics described below.

 Without going into further detail, we present \red{two statements (of Borsuk-Ulam
 type): Assertion \ref{conj:far} and Conjecture \ref{conj:close}. Together,
 they imply Conjecture \ref{conj:main}, as asserted in
 \cite[p.~40]{simon03:_habil}, compare \cite{SL}. This follows rather easily
 using the methods of \cite{SST2}}. To state them, we first make a definition.

 \begin{definition}
   \begin{enumerate}
\item  If $A$ is a subset of an affine space,
     $\co(A)$ stands for its convex hull.
   \item If $F\subset X\times Y$ and $y\in Y$, then $F^{-1}(y):=\{x\in X\mid
     (x,y)\in F\}$.
   \item If $F\subset \Delta(L) \times Y$ we define
     $cF\subset\Delta (L) \times Y$ by $cF:= \{ (x,y)\ |\ x\in {\co} (F^{-1}(y)
     )\}$.
   \item If $F\subset \Delta (L)\times Y$  and $Y\subset
     {\reals} ^L$ then is $F ^+\subset\Delta
     (L)\times Y$ defined by
     \begin{equation*}
F^+:=\{ (p,y)\mid \exists
     (p,x)\in F:\, x^l\leq y^l\, \forall l\in L \text{ and } x^l=y^l\text{ if }
     p^l>0\}.
   \end{equation*}
We call $F^+$ the \emph{$Y$- saturation} of $F$.  If $F=F^+$ then $F$ is
called saturated.
\item A constant correspondence $\Delta (K)\ni p\mapsto U$
is denoted $U^{cst}$ (that is, $U^{cst}:=\Delta (K)\times U$).
   \end{enumerate}
  \end{definition}

\def\L{{\cal L}}

In both the Assertion \ref{conj:far} and Conjecture \ref{conj:close} below we
assume 
that $I=[a,b]\subset \reals $ is a non-trivial closed segment,
$ K$  is a finite set and $\L $ is a family of non-empty proper subsets of $K$
which covers $K$. Each simplex $\Delta (L),\,L\in \L$,  is considered
as a subset of $\Delta (K)$.

\begin{assertion}\label{conj:far}
Suppose for every
$L\in \L\cup \{K\} $ there is given a saturated correspondence
$F_L\subset \Delta (L)\times I ^L$  with property ${\M} $ for $\Delta (L)$ and  a
closed convex set $U_L\subset I ^K$ containing the point
$v_+=(b,...,b)$.  Set 
\begin{equation}\label{eq:FLtilde}
\widetilde{F _L}=
 F_L\times I ^{K\setminus L}\subset \Delta (K)\times I ^K \mbox{ for } L\in \L
\end{equation}
 Also let $U:=\bigcap_{L\in \L \cup \{K\}}U_L, \ F:=U^{cst}\cap F_K \mbox{ and }
 G:=U^{cst}\cap \bigcup_{L\in \L} \widetilde{F_L}$, and
define a correspondence $\Gamma\subset \Delta (K)\times I ^K$
so that for $y\in I^K$
$$\Gamma ^{-1}(y)=\co (G^{-1}(y))\cup \bigcup _{x\in F^{-1}(y)}\co (\{x\}\cup
G^{-1}(y)) 
$$
The assertion is that if $\widetilde{F _L} \subset U_L^{cst}$ for each $L\in \L\cup \{K\}$,
then $\Gamma $ has property $\M$ for $\Delta (K)$.
\end{assertion}

The key property in the definition of $\Gamma$ is that for every $y\in
\Delta(K)$
we convexify sets in
$F^{-1} (y) \cup G^{-1} (y)$ that have at most one member $x$ whose projection
to $\Delta(K)$ belongs to $\Delta (K) \setminus  \bigcup_{L\in \L}\Delta (L)$.

 \begin{conjecture}\label{conj:close}
This time suppose the saturated correspondences $F_L\subset \Delta (L)\times I
^L$ 
with property ${\M} $ for $\Delta (L)$ are given merely for those sets $L\in \L $
which are maximal with respect to inclusion. For  all $L\in \L $ let
$$F_L:=\bigcup_{L\subset J\in \max \L}   F_J \cap
   (\Delta (L)\times I ^L)\subset \Delta(L)\times I ^L$$ and define
   $\tilde {F} _L$ as in \eqref{eq:FLtilde} and
   $\Gamma\subset \Delta (K)\times I ^K$ so that for $y\in I^K$
  $$\Gamma^{-1}(y):= \bigcup
\{\co(x_{L_1}, \dots , x_{L_m}) \mid
{x_{L_i} \in \tilde F_{L_i}^{-1}(y),\ L_i = L_j \in \max \L \Rightarrow \
     i=j}\} .$$
     The conjecture is  that if $L_1\cap L_2\in \L$ for all $L_1, L_2\in \L$, then
     $\Gamma $
   has property $\M$ for $\Delta(K)$.
 \end{conjecture}
  Because ---in a way similar to Assertion \ref{conj:far}--- in Conjecture
  \ref{conj:close} we have a
  restriction on how we convexify sets to obtain $\Gamma$, its character is close to the original 
  Borsuk-Ulam Theorem (where one convexifies using opposite points).  \red{It suggests the statement of
%
%
%
  Corollary \ref{theo:special_fam_bors}, proved above. Indeed, this relation
to Conjecture \ref{conj:close} was one of the main motivations for the work
presented in this paper.}

\red{We believe that we can prove Assertion \ref{conj:far}. We plan to publish
  this in a forthcoming paper. Conjecture \ref{conj:close} is still
  open. Provided we can establish it, as well, we will eventually put all the
  details together and present a proof of the game theoretic Conjecture \ref{conj:main}.}

\bibliographystyle{plain}

\end{document}